\documentclass{amsart}
\usepackage{indentfirst,latexsym,bm,amsmath,amssymb,amsfonts,lscape,rawfonts,eufrak}
\usepackage{algorithm,algorithmic}
\usepackage{graphicx}

\newcommand{\R}{\mathbb{R}}

\newtheorem{theorem}{Theorem}[section]
\newtheorem{proposition}[theorem]{Proposition}

\newtheorem{definition}[theorem]{Definition}

\newtheorem{Rk}[theorem]{Remark}

\numberwithin{equation}{section}

\begin{document}

\title{A New Multiscale Representation for Shapes and Its Application to Blood Vessel Recovery}

\author{Bin Dong$^1$ \and Aichi Chien$^2$ \and Zuowei Shen$^3$  \and Stanley
Osher$^1$}

\date{}

\thanks{
1. Department of Mathematics, University of California, Los
Angeles, CA, (supported by NIH GRANT, P20 MH65166; NSF GRANT,
DMS-0714807);}

\thanks{
2. Division of Interventional
Neuroradiology, David Geffen School of Medicine at UCLA, 10833
LeConte Ave, Los Angeles, CA, (supported by NSF GRANT,
CCF-0830554);}

\thanks{
3. Department of Mathematics, National
University of Singapore, 2 Science Drive 2, Singapore, 117543.}

\maketitle

\begin{abstract}
In this paper, we will first introduce a novel multiscale
representation (MSR) for shapes. Based on the MSR, we will then
design a surface inpainting algorithm to recover 3D geometry of
blood vessels. Because of the nature of irregular morphology in
vessels and organs, both phantom and real inpainting scenarios
were tested using our new algorithm. Successful vessel recoveries
are demonstrated with numerical estimation of the degree of
arteriosclerosis and vessel occlusion.
\end{abstract}

\section{Introduction}

\subsection{Literature Reviews and Motivations}

Multiscale representation (MSR) of functions, e.g. wavelets, has
been extensively studied in the past twenty years \cite{Dau,
Mallat}. However, when one deals with shapes, e.g. biological
shapes in $\R^3$, most of the classical theories and algorithms
cannot be directly extended. In this paper, we will propose a new
MSR for shapes based on PDEs and level set method. Although we
shall focus on studying 3D biological shapes/surfaces, the MSR
that we introduce here applies to general shapes/surfaces in both
2D and 3D.

Many attempts have been made in the past on designing
wavelet-typed MSR for 3D shapes \cite{NHBT,NHBT1,Grshields,BDong}.
Among them, the method proposed by Nain et. al. \cite{NHBT,NHBT1}
is especially effective to study biological shapes. They first map
the shape (triangulated) onto the unit sphere so that one obtains
a vector-valued function $f: \mathbb{S}^2\mapsto\R^3$; then apply
spherical wavelet decomposition \cite{SS} to each component of
$f$. However, the wavelet coefficients are not intrinsic to the
shape, but dependent on the mapping $f$. Furthermore, finding a
good mapping from a shape to the unit sphere (or to some other
canonical domains) is nontrivial and in fact a popular ongoing
research area (see e.g. \cite{GWCTY, WGCTY, GWCTY1, JLG, HATKSH,
GGS, GY, JWY, STDOT, PH}).

Another interesting approach was proposed by Pauly et. al.
\cite{PKG}, where they introduced an MSR for point-based surfaces.
Their idea was to use Moving Least Square method \cite{Levin} to
define a series of smoother and smoother point-based surfaces, and
then define wavelet coefficients as the displacements from two
successive levels. Their method only requires a local
parametrization of the point-based surface which is easy to
calculate. However, the application of their method is rather
limited in medical image analysis, because most of the biological
shapes are not point-based.

Motivated by Pauly et. al.'s work, we will propose a new MSR for
shapes in Section \ref{Section:LSBMSR}. The basic idea is using
level set motions via solving some properly chosen Hamilton-Jacobi
(HJ) like equation to obtain a sequence of shapes that become
smoother and smoother as time evolves (analogous to coarse level
approximation in wavelet decomposition). Then we carefully define
the so-called ``details" (analogous to wavelet coefficients) of
the MSR which carry important geometric information and facilitate
a perfect reconstruction. While the wavelet based multi scale
decomposition and reconstruction use filters, which are linear
processes, the proposed new MSR for shapes uses (nonlinear) PDEs
for both decomposition and reconstruction.  However, the spirit is
the same, i.e. separate features from smooth components of the
surface. Due to the level set formulation, parametrization is no
longer needed.

\subsection{Shape Modelling and Evolution PDEs}

Throughout this paper, shapes are defined to be smooth boundaries
of domains $\Omega\in\R^3$ and are represented by level set
functions, typically signed distance functions. We note, however,
that point-based and triangulated surfaces can also be handled in
a similar way, as noted in (4) of Remark \ref{Remark:Def:MSR}.

A level set function $\phi$ that represents the shape
$\partial\Omega$ is defined as follows
\begin{equation*}
\phi(x)\left\{
\begin{array}{ll}
<0 & x\in\Omega;\\
>0 & x\in\Omega^{c}.
\end{array}
\right.
\end{equation*}
We always assume that the function $\phi$ is at least Lipschitz
continuous.

Level set motions can be achieved by solving the following HJ like
equation \cite{OS},
\begin{equation}\label{Eqn:LSMotion}
\phi_t+v_n\left(\nabla\phi\right)|\nabla\phi|=0,\quad
\phi(x,0)=\phi_0(x),
\end{equation}
where we take $(x,t)\in\mathcal{D}\times[0,T]$ with $\mathcal{D}$
some bounded domain in $\R^3$ and $T>0$. Here
$v_n\left(\nabla\phi\right)$ is the normal velocity, which
essentially depends on $\nabla\phi$ while second order derivatives
of $\phi$ may be involved (e.g. mean curvature). If it only
depends on first order derivatives, then it is a HJ equation. We
also assume that the PDE \eqref{Eqn:LSMotion} is geometric
\cite{CGG,Giga}, which guarantees contrasts invariance.
Comprehensive theoretical analysis of PDE \eqref{Eqn:LSMotion} and
surface evolution equations can be found in
\cite{CGG,Giga,ES1,ES2,ES3,ES4,CL,CL1,CIL}.

The choice of velocity fields is very important and yet very
non-unique. We need to choose one that generates a ``meaningful"
MSR , e.g. a sparse representation, for a given smooth shape.
Generally speaking, we want the zero level set of $u(x,t)$ becomes
smoother and smoother as $t$ increases. This is in fact a typical
scale space behavior that has been studied for decades (see e.g.
\cite{AGLM, Morel}). It is known \cite{AGLM, Morel} that under
some general axiomatic hypothesis and some invariance (i.e.
rotation and contrast invariance) assumptions on
$\{u(x,t)\}_{t\ge0}$, $u(x,t)$ must be a viscosity solution to a
PDE of the form \eqref{Eqn:LSMotion}, with the velocity field
$v_n$ only depending on the principle curvatures of level sets of
$u$ and time $t$. In other words, a ``meaningful" velocity field
must be curvature dependent.

The velocity fields that we shall focus in this paper are
\begin{equation}\label{Eqn:Velocities}
v_n=c-\lambda\kappa,\ \lambda>0\quad\mbox{and}\quad
v_n=\kappa_a-\kappa,
\end{equation}
where $c\in\R$ is some constant, $\kappa$ is the mean curvature
defined as $\kappa:=\nabla\cdot\frac{\nabla\phi}{|\nabla\phi|}$,
and $\kappa_a$ is the average mean curvature \cite{ESvpmc}. Note
that for $v_n=\kappa_a-\kappa$, the PDE \eqref{Eqn:LSMotion}
generates an volume preserving mean curvature motion
\cite{ESvpmc,BS,RW,RS}.

\section{Level Set Based MSR of Shapes: Continuous Transforms and Discrete
Algorithms}\label{Section:LSBMSR}

Let $\Omega_t\in\R^3$ be some domain with scale $t$, and
$S_t:=\partial\Omega_t$ be the shape at scale $t$ represented by
some time-dependent level set function $\phi(x,t)$, i.e.
$\phi(x,t)<0$ for $x\in\Omega_t$, $\phi(x,t)>0$ for
$x\in\Omega_t^{c}$, and
\begin{equation}\label{Eqn:Def:St}
S_t=\{x\in\R^3\ |\ \phi(x,t)=0\}_{t\ge0}.
\end{equation}
Here $S_0$ denotes the original shape with the corresponding level
set function $\phi_0(x)=\phi(x,0)$. Throughout the rest of the
paper, the function $\phi(x,t)$ is always taken to be the solution
of \eqref{Eqn:LSMotion}. For some properly chosen $v_n$ in
\eqref{Eqn:LSMotion}, e.g. with $v_n=-\kappa$ or
$\kappa_a-\kappa$, we can obtain a continuous series of shapes
$\{S_t\}_{t\in[0,T]}$, which tends to become smoother when $t$
increases. Based on this, we define our continuous level set based
MSR of $S_0$ as follows.

\begin{definition}\label{Def:CMR:Alt}
Let $\phi(x,t)$ be the solution of the PDE \eqref{Eqn:LSMotion}
and $(x,t)\in\mathcal{D}\times[0,T]$. We now understand $x_l(t)$
as a path on the propagating $l$-th level set of $\phi$, i.e.
$\phi(x_l(t),t)=l$. For simplicity, we shall omit the subscript
``$l$" unless a particular level set is considered.

\begin{enumerate}

\item We now define the \textbf{multiscale transformation (MST)}
of $\phi_0(x)$ as
\begin{equation}\label{Eqn:CMR:Decomp:LSB}
\vec
W(x,t):=W(\phi_0):=-v_n\frac{\nabla\phi}{|\nabla\phi|}=-x'(t).
\end{equation}
Vector $-x'(t)$ is the \textbf{displacement vector} and
$w(x,t):=-v_n(x,t)$ is the \textbf{detail} of the MST.

\item We shall call $\vec W(x,t)$ the \textbf{displacement vector
field} at scale $t$, and denote $\vec W_|(x,t)$ ($w_|(x,t)$) as
the restriction of $\vec W(x,t)$ ($w(x,t)$) on $S_t$.

\item The MSR for the original shape $S_0$ in terms of $\phi_0(x)$
is denoted as
$$\mbox{MSR}(\phi_0)=\Big\{\{\vec W(x,t)\}_{t\in(0,T)},\phi(x,T)\Big\}.$$

\item We define the \textbf{inverse multiscale transformation
(IMST)} via solving the following PDE
\begin{equation}\label{Eqn:CMR:Recons:LSB}
\psi_\tau+\vec W(x,T-\tau)\cdot\nabla\psi=0,\quad
\psi(x,0)=\phi(x,T).
\end{equation}
for given $T>0$ and $0\le\tau\le T$.
\end{enumerate}
\end{definition}

\begin{Rk}\label{Remark:Def:MSR}

\hspace*{3in}
\begin{enumerate}

\item The technique of generating a sequence of the spaces
$\{S_t\}$ via solving PDEs is known as scale space decomposition
(see e.g. \cite{AGLM, Morel}). However, a classical scale space
analysis does not study the details as defined in (2), and does
not have a reconstruction as in (4).

\item The last identity in \eqref{Eqn:CMR:Decomp:LSB} can be
easily shown by using PDE \eqref{Eqn:LSMotion} and the assumption
that $x'(t)$ is aligned with normal directions of level sets of
$\phi$.

\item The detail $w_|(x,t)$ is a function on $S_t$ that
characterizes intrinsic geometric information of the shape at
scale $t$. Here by intrinsic we mean that $w_|(x,t)$, as well as
$\{S_t\}_{t>0}$, does not depend on the initial embedding $\phi_0$
for a large class of functions \cite{Giga}, but only depends on
$S_0$. Therefore, we now have an intrinsic MSR for $S_0$:
\begin{equation}\label{Eqn:MSR}
\mbox{MSR}(S_0)=\Big\{\{\vec W_|(x,t)\}_{t\in(0,T)},S_T\Big\}.
\end{equation}
Furthermore, the above MSR is invariant under translation and
rotation of $S_0$.

\item The MSR defined above can be easily adapted to a point-based
or triangulated surface. One simply need to first associate the
surface with a level set function and then perform the MST. For
point-based surfaces, the IMST from its MSR \eqref{Eqn:MSR} can be
point-wise defined as $S_0=S_T+\int_0^T \vec W_|(x,t)dt$ or
equivalently $x_0(0)=x_0(T)+\int_0^T -x_0'(t)dt$, which is
obviously true.

\end{enumerate}
\end{Rk}

Now the question is that if we have perfect reconstructions via
\eqref{Eqn:CMR:Recons:LSB}. The answer is given in the following
proposition, which directly follows from theories of ODEs.

\begin{proposition}\label{Prop:LSBMSR}
Assume that $\vec W(x,t)$ stays Lipschitz continuous for
$(x,t)\in\mathcal{D}\times[0, T]$. Then the equation
\eqref{Eqn:CMR:Recons:LSB} \textbf{inverts} the MST defined by
\eqref{Eqn:CMR:Decomp:LSB} in the sense that
$\psi(x,\tau):=\phi(x,T-\tau)$ is the unique solution of
\eqref{Eqn:CMR:Recons:LSB}.
\end{proposition}

\begin{Rk}\label{Remark:Prop}
\hspace*{3in}
\begin{enumerate}

\item The assumption in Proposition \ref{Prop:LSBMSR} is not
always valid (e.g. $v_n=c<0$ and $\phi_0(x)$ representing a cube).
However, if we choose, for example, $v_n=-\kappa$,
$v_n=\kappa_a-\kappa$ or $v_n=c-\lambda\kappa$ (for $\lambda>0$),
and choose some appropriate ending time $T>0$ (e.g. before any
topological changes occur), the above assumption will be valid and
we will have a perfect reconstruction using
\eqref{Eqn:CMR:Recons:LSB} \cite{Giga,ES}.

\item Generally speaking, the vector field $\vec W(x,t)$ does not
stay Lipschitz globally in time. This happens precisely when the
corresponding surface evolution starts to develop singularities,
and it is difficult to find a surface evolution that guarantees to
have global smooth solutions for a general initial surface $S_0$.
For some special class of initial surfaces, however, it is
relatively easy to find such motion. Taking $v_n=\kappa_a-\kappa$
for example, it is shown in \cite{ESvpmc} that if the initial
surface $S_0$ is close enough to a certain sphere (but not
necessarily convex), then $S_t$ stay smooth and eventually
converges exponentially fast to a sphere with the same volume as
$S_0$.
\end{enumerate}
\end{Rk}

Notice from Definition \ref{Def:CMR:Alt} and Proposition
\ref{Prop:LSBMSR} that to perfectly reconstruct $\phi_0(x)$ from
$\phi(x,T)$, we need to store the entire vector field $\vec
W(x,t)$ for every $x\in\mathcal{D}$ and all scale $t$. However, in
practice, we only want a perfect reconstruction of $S_0$, and thus
we do not need that much information. Therefore, only the
displacement vectors within a narrow band of the zero level set of
$\phi(x,t)$ need to be stored.

We can be even more ``greedy" here by only storing $\vec
W_|(x,t)$. When performing inverse transform, we will need to
extend $\vec W_|(x,t)$ to at least a narrow band of the zero level
set of $\phi(x,t)$. Note that no extension can guarantee an exact
recovery of the vector field $\vec W(x,t)$, and hence the
reconstruction of $S_0$ will not be exact. However, if the
extension is conducted accurately and the mesh grid is dense
enough, i.e. the resolution of the shape is high enough, the
reconstruction should be more and more accurate. The extension we
shall adopt here is such that the extended vectors are constant in
the normal directions of each level set of $\phi(x,t)$
\cite{BCSO}. For simplicity, we will use a local search method to
extend $\vec W_|(x,t)$ to a narrow band of the zero level set of
$\phi(x,t)$.

Our proposed discrete version of MSR is given in Algorithm
\ref{MSTIMS LSB_Alg}.

\begin{algorithm}

\caption{Level Set Based MST and IMST \label{MSTIMS LSB_Alg}}
\begin{algorithmic}
\STATE Start from the given level set function $\phi_0(x)$
representing shape $S_0$. Choose time steps
$0=t_0<t_1<\ldots<t_N=T$, where $\max_i(t_{i+1}-t_i)$ is small.

\textbf{Initialize:} Sample a point set $X_0$ from $S_0$ (either
uniformly or non-uniformly).

{\hspace*{2in}\textbf{MST:}}

\WHILE{$i\le N$}

\STATE{\textbf{1.}} Starting from $\phi(x,t_{i-1})$, solve PDE
\eqref{Eqn:LSMotion} for $t\in[t_{i-1},t_i]$ and obtain
$\phi(x,t_i)$.

\STATE{\textbf{2.}} Orthogonally project $X_{i-1}$ onto the zero
level set of $\phi(x,t_i)$ and obtain $X_i$.

\STATE{\textbf{2.}} Compute the discrete displacement vector by
$\vec W_{|i}=X_i-X_{i-1}$, and $i\leftarrow i+1$.

\ENDWHILE

We then obtain the discrete MSR of $S_0$: $\mbox{MSR}(S_0):=\{\vec
W_{|1}, \vec W_{|2}, \ldots,\vec W_{|N}, \phi(x,T)\}$.

{\hspace*{2in}\textbf{IMST:}}

\STATE{\textbf{1.}} Extend the vector fields $\{\vec
W_{|i}\}_{i=1}^N$ such that the values are constant along normal
directions of the level sets of $\phi(x,t_i)$.

\STATE{\textbf{2.}} Solve \eqref{Eqn:CMR:Recons:LSB} using $\vec
W_{|i}$ within interval $[t_i, t_{i-1}]$ iteratively for each $i$.

\end{algorithmic}

\end{algorithm}

\section{Numerical Experiments on the
MSR}\label{Section:LSBMSR:Numerical}

One of the key steps of implementing Algorithm \ref{MSTIMS
LSB_Alg} is to solve the evolution PDE \eqref{Eqn:LSMotion}
efficiently. There are many ways of solving equation
\eqref{Eqn:LSMotion}. The most straightforward way is to use
monotone finite difference schemes \cite{OS, OF}. However, it is
not very efficient computationally. To overcome this, Merriman,
Bence and Osher introduced a diffusion-based level set motion in
\cite{MBO1, MBO2}, and it was further studied in \cite{Ishii, R,
RM, RMO}, where in \cite{Ishii} the correctness of the method is
rigorously proven. In \cite{RW}, Ruuth and Wetton introduced a
fast algorithm to calculate volume preserving motion by mean
curvatures. All these methods speeded up curvature driven motions
drastically.

In this section, we will recall the fast algorithms of level set
motion for the case $v_n=c$ and $v_n=\kappa_a-\kappa$ given by
\cite{RW, MBO1, MBO2, RM}. These algorithms will be used in the
later sections to generate fast multiscle decompositions of
shapes.

Now we first recall the fast method of solving
\eqref{Eqn:LSMotion} with $v_n=c$ (see \cite{MBO1, MBO2, RM}) in
Algorithm \ref{MotionInNormalDirection_Alg}.

\begin{algorithm}
\caption{Level Set Motion with Constant Normal Velocity
\label{MotionInNormalDirection_Alg}}
\begin{algorithmic}
\STATE Start from a given shape represented by $\phi$.

\WHILE{$t<T$}

\STATE{\textbf{1.}} Define the corresponding characteristic
function by $\chi=\textbf{1}_{\{\phi<0\}}$. Set $V_0$ equal to the
volume of $\{\phi<0\}$.

\STATE{\textbf{2.}} Starting from $\chi$, evolve $\bar\chi$ for a
time $\Delta t$ by $\bar\chi_t=\nabla^2 \bar\chi$.

\STATE{\textbf{4.}} Sharpen:
\begin{equation*}
\chi=\left\{
\begin{array}{ll}
1& \mbox{if}\ \bar\chi>0\\
0& \mbox{otherwise}
\end{array}
\right.
\end{equation*}

\STATE{\textbf{5.}} Let $t\leftarrow t+\Delta t$. Compute
$\phi(x,t)$ from $\chi$ via fast sweeping method \cite{TCOZ}.

\ENDWHILE
\end{algorithmic}
\end{algorithm}

We now recall the fast implementation of \eqref{Eqn:LSMotion} with
$v_n=\kappa_a-\kappa$ proposed by Ruuth and Wetton \cite{RW} in
Algorithm \ref{VPMeanCurvMotion_Alg}. Their algorithm is based on
the diffusion-based mean curvature motion proposed by \cite{MBO1,
MBO2}. Note that if we remove step 3 in Algorithm
\ref{VPMeanCurvMotion_Alg} and choose $\lambda=0.5$ in step 4, it
is exact the fast mean curvature motion proposed in \cite{MBO1,
MBO2}.

\begin{algorithm}
\caption{Volume Preserving Mean Curvature Motion:
$v_n=\kappa_a-\kappa$. \label{VPMeanCurvMotion_Alg}}
\begin{algorithmic}
\STATE Start from a given shape represented by $\phi$.

\WHILE{$t<T$}

\STATE{\textbf{1.}} Define the corresponding characteristic
function by $\chi=\textbf{1}_{\{\phi<0\}}$. Set $V_0$ equal to the
volume of $\{\phi<0\}$.

\STATE{\textbf{2.}} Starting from $\chi$, evolve $\bar\chi$ for a
time $\Delta t$ by $\bar\chi_t=\nabla^2 \bar\chi$.

\STATE{\textbf{3.}} Determine the threshold value that preserves
the volume of the set: i.e. find a $0<\lambda<1$ s.t.
$$\Big||\{x:\bar\chi<\lambda\}|-V_0\Big|<\varepsilon.$$

\STATE{\textbf{4.}} Sharpen:
\begin{equation*}
\chi=\left\{
\begin{array}{ll}
1& \mbox{if}\ \bar\chi>\lambda\\
0& \mbox{otherwise}
\end{array}
\right.
\end{equation*}

\STATE{\textbf{5.}} Let $t\leftarrow t+\Delta t$. Compute
$\phi(x,t)$ from $\chi$ via fast sweeping method \cite{TCOZ}.

\ENDWHILE
\end{algorithmic}
\end{algorithm}

\begin{figure}[htp]
\centering
    \includegraphics[width=5.0in]{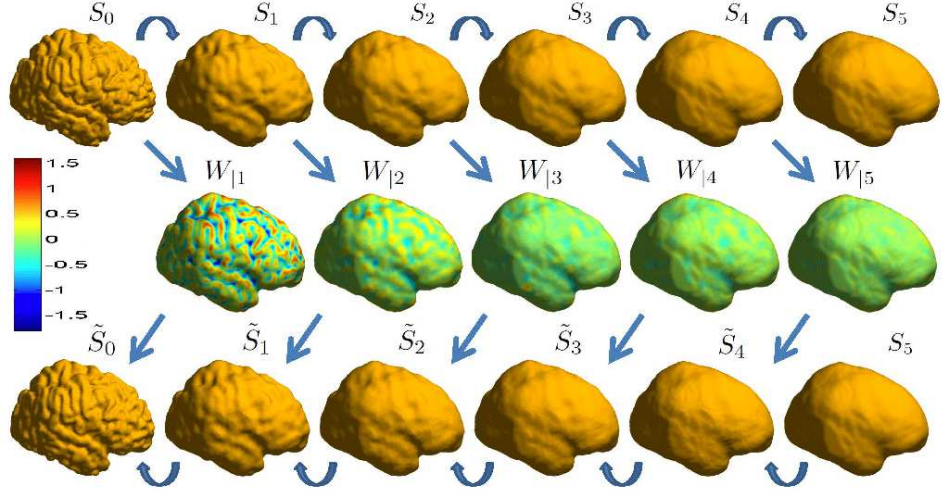}
\caption{First row (left to right): MST $S_0, S_1, \ldots, S_5$.
Second row shows the details of MSR on $S_1, \ldots, S_5$. Third
row shows IMST $\tilde S_i$, $i=0,1,\ldots,4$, where the Hausdorff
distance between $S_i$ and $\tilde S_i$ are: $1.12h$, $0.74h$,
$0.74h$, $0.69h$, and $0.63h$ respectively (with $h$ the mesh
size).}\label{Fig:MST:IMST:cortex}
\end{figure}

\begin{figure}[htp]
\centering
    \includegraphics[width=0.9in]{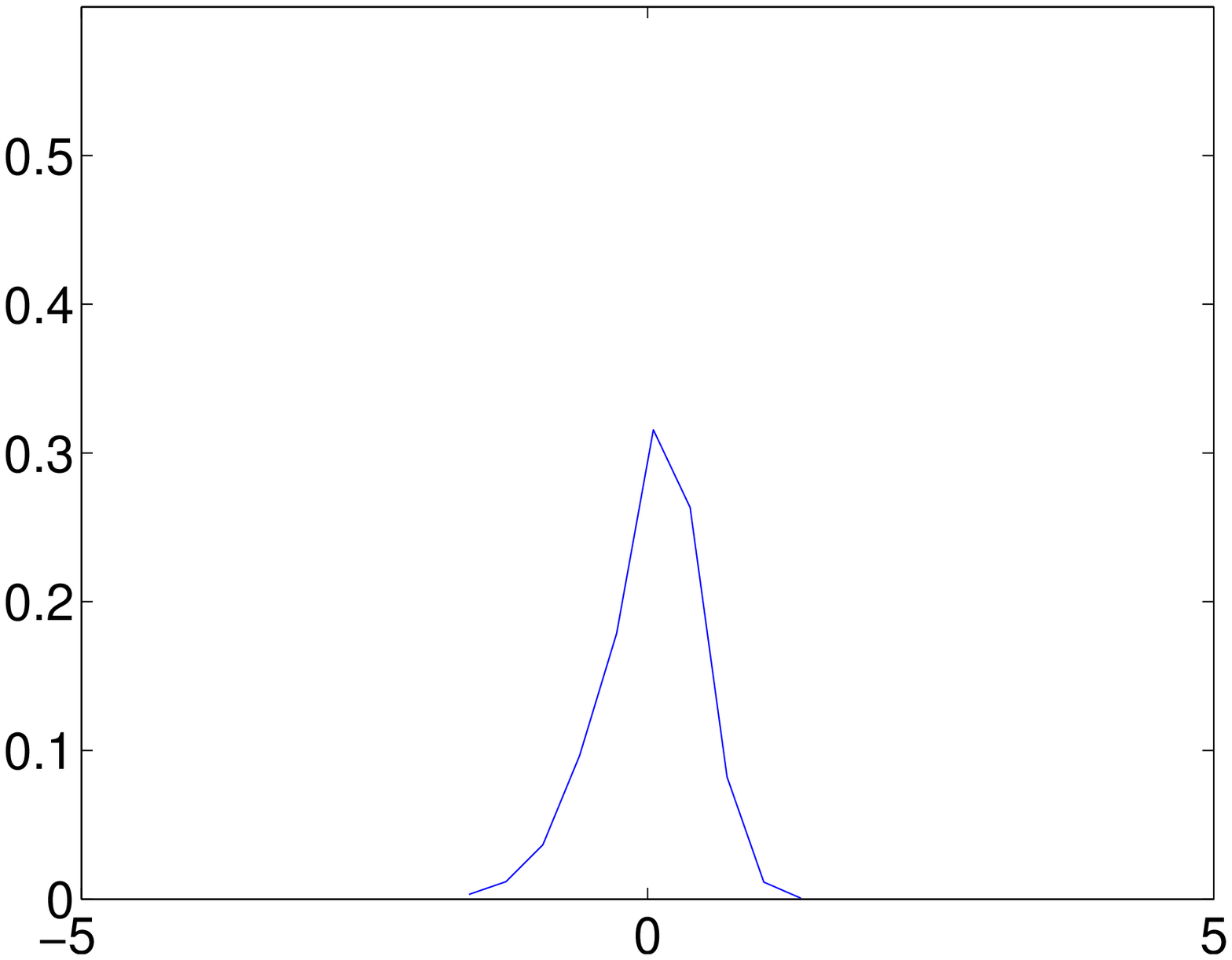}
    \includegraphics[width=0.9in]{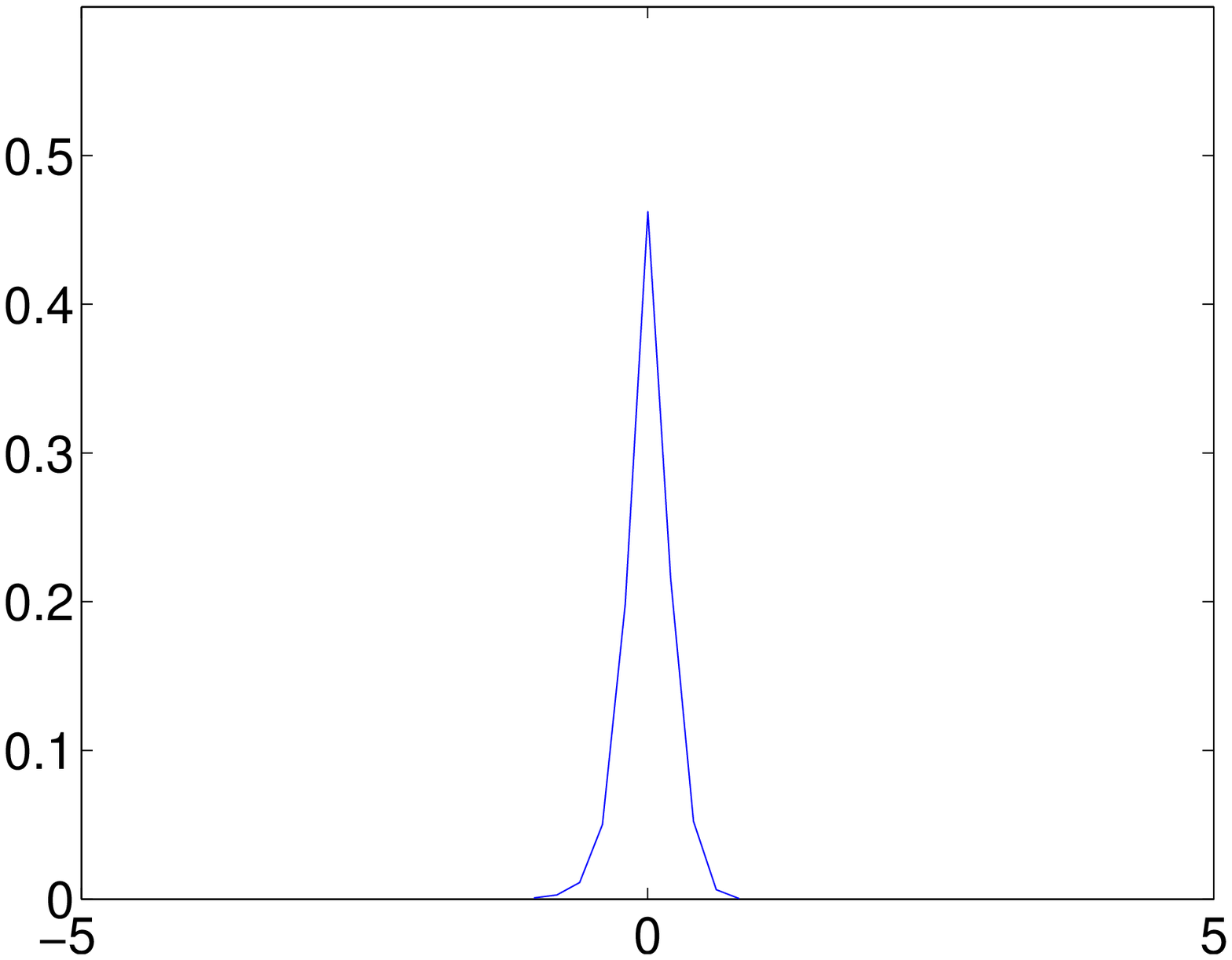}
    \includegraphics[width=0.9in]{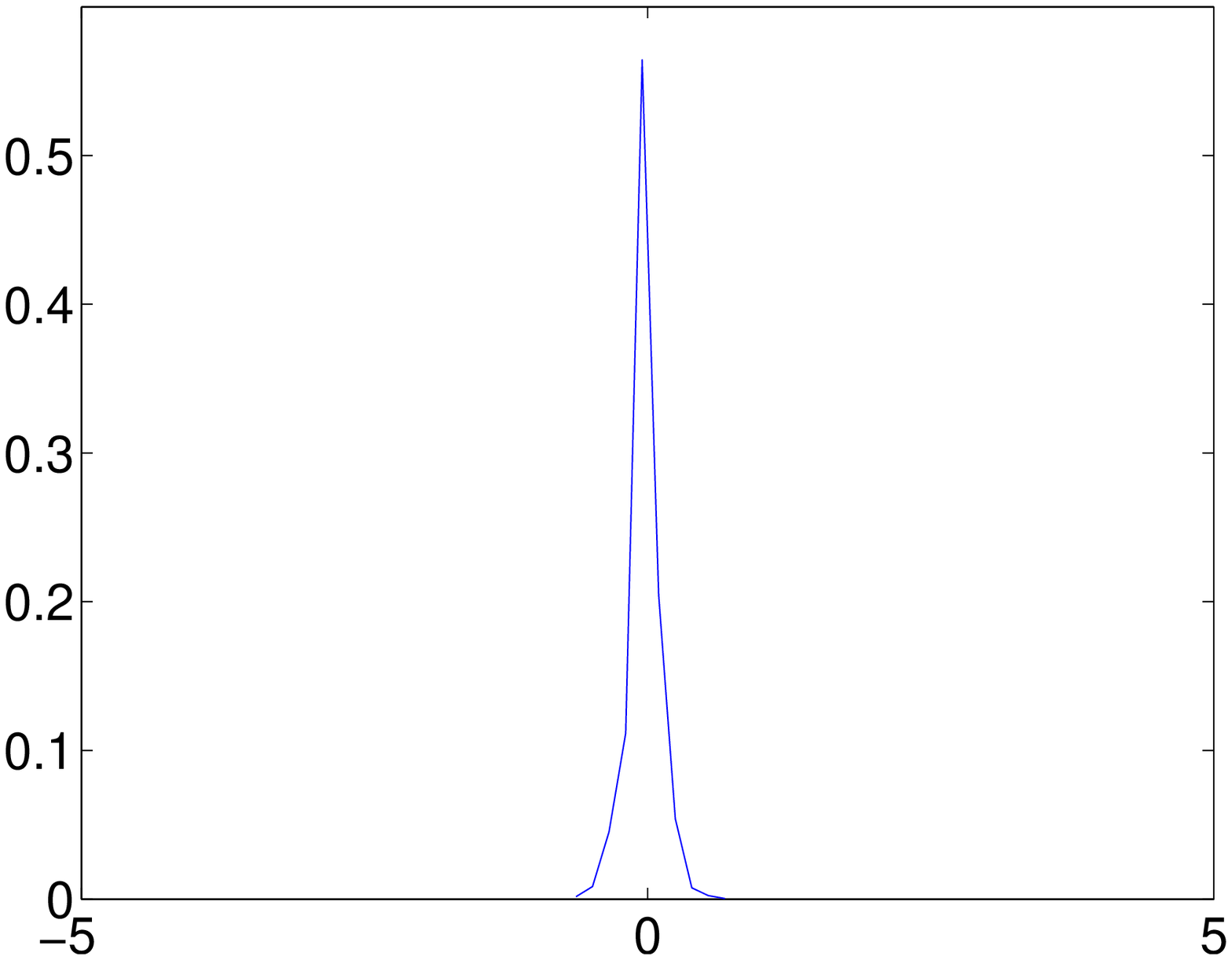}
    \includegraphics[width=0.9in]{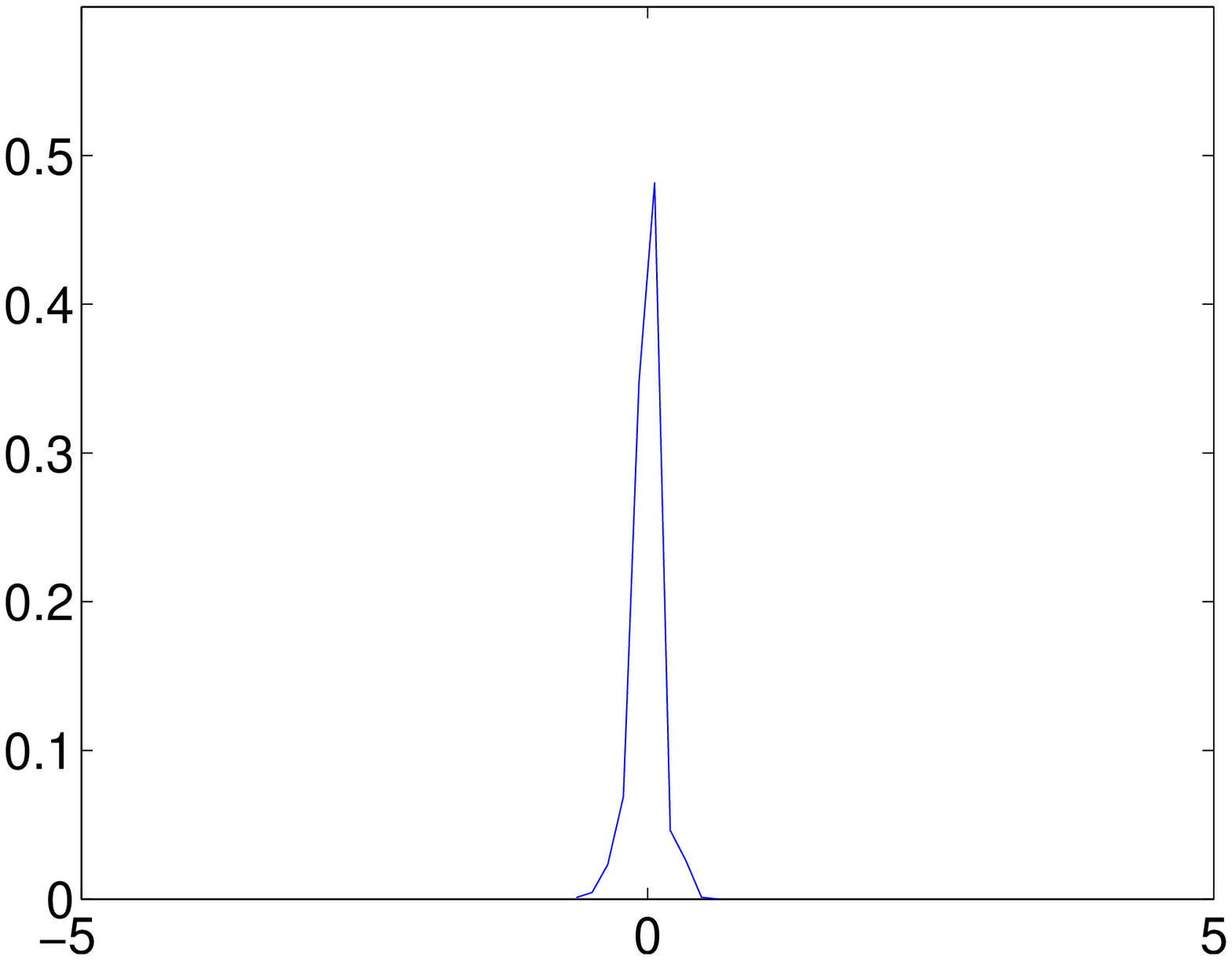}
    \includegraphics[width=0.9in]{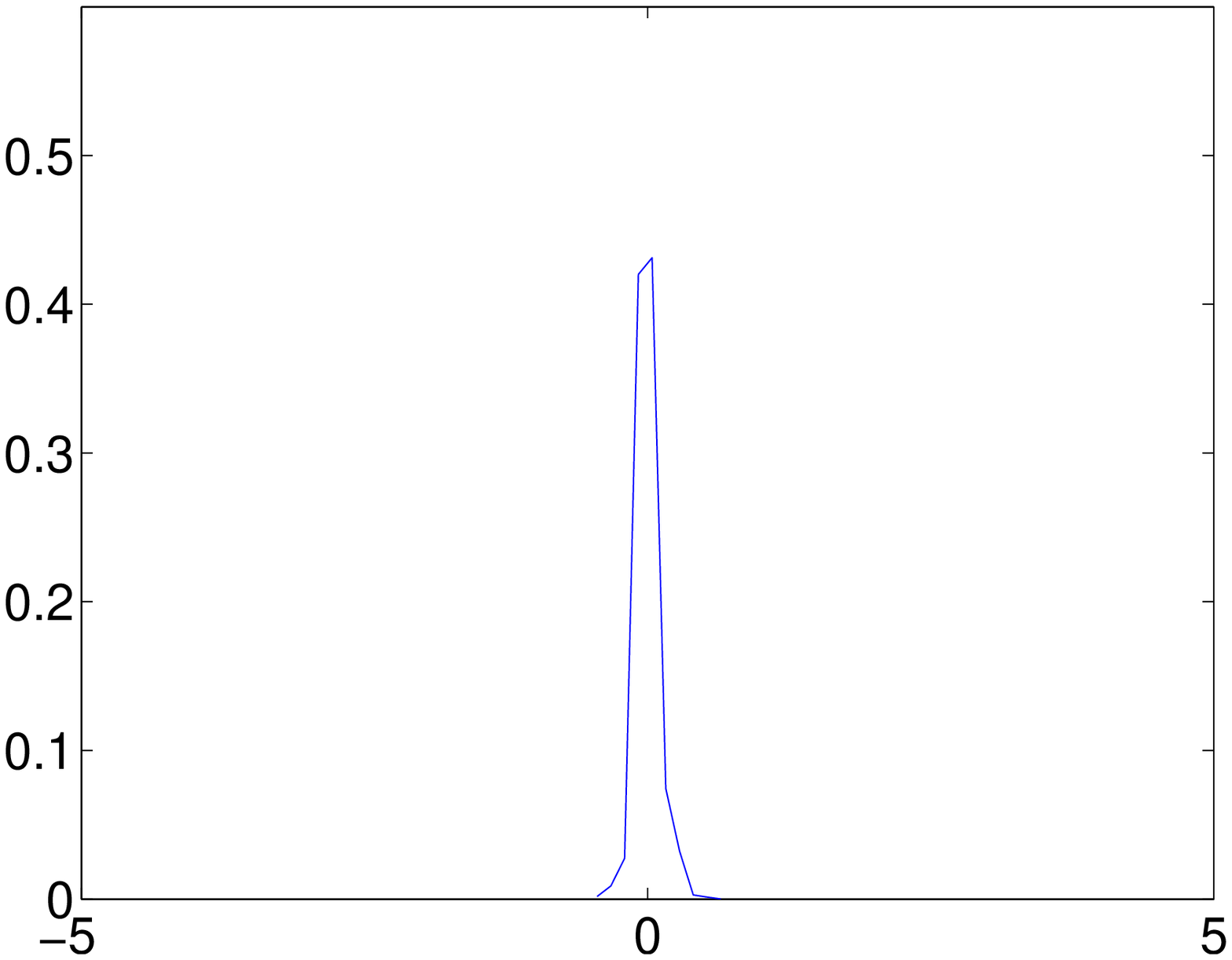}
\caption{Histograms of $\vec W_{|i}$ for $i=1,\ldots,5$ (left to
right).}\label{Fig:Sparseness}
\end{figure}

Some numerical results of the MST and IMST in Algorithm
\ref{MSTIMS LSB_Alg} are presented in Figure
\ref{Fig:MST:IMST:cortex} using the tested shape (right hemisphere
of a cortex). The velocity field is chosen to be
$v_n=\kappa_a-\kappa$ and 5 levels of decomposition are conducted
(first and second row of Figure \ref{Fig:MST:IMST:cortex}).
Details $\vec W_{|i}$ are drawn on the surface $S_i$ (second row
of Figure \ref{Fig:MST:IMST:cortex}), where the value is positive,
when $\vec W_{|i}$ is pointing outwards and negative when it is
pointing inwards. The IMST is also presented in Figure
\ref{Fig:MST:IMST:cortex} where $\tilde S_i$ denotes the
reconstruction of level $i$ from level $i+1$. As we can see,
although the reconstructions are not exact for each level, they
are quite accurate in the sense that most of the features are well
reconstructed. We also illustrate sparseness of the coefficients
$\{\vec W_{|i}\}_{i=1}^5$ in Figure \ref{Fig:Sparseness}, where
one can see that the energy of $\vec W_{|i}$ are concentrated
around 0, especially for the later levels.

\section{Application in Blood Vessel Recovery}\label{Section:MSRInpainting}

Evaluating missing parts in medical images provides important
information as signs of diseases. One of the most common situation
is the phenomenon of vessel narrowing or occlusion in angiographic
images. Estimating these abnormalities can help document disease
progression.

The recovery of image surfaces such as blood vessels can be
regarded as a surface inpainting problem \cite{DMGL,VCBS}.
Inpainting problems, for both images and surfaces, have been
extensively studied in the literature \cite{NRC, CKSip, BSCB, BBS,
CSip, BEG1, BEG2, ESip, DBip, BVSO, CCSS, CCSip, ESQD, CSZ,DMGL,
VCBS, ABK, BMRST, Whitaker, WSI,ZOMK}. They occur when part of the
data in an image or regions of a surface is missing or corrupted.
The major task of inpainting is to fill in the missing information
based on the geometry of the image/surface. In this section, we
will propose a new surface inpainting algorithm based on the MSR
in Section \ref{Section:LSBMSR} for blood vessel reconstruction
that arises in medical image analysis.

Our surface inpainting algorithm (Algorithm
\ref{Inpainting_Alg_LSB} below) inherits the following
framelet-based image inpainting structure proposed by Cai et. al.
\cite{CCSip}:
\begin{enumerate}
\item Take framelet transform of the given image;

\item Truncate the framelet coefficients via soft-thresholding and
reconstruct;

\item Apply the exact data outside the inpainting domains, and
repeat.
\end{enumerate}
Since we already have an MSR for surfaces, the first step above
can be replace by our MST. For the second step, we shall solve the
following PDE for IMST instead of the PDE
\eqref{Eqn:CMR:Recons:LSB} that was originally proposed in
Definition \ref{Def:CMR:Alt}:
\begin{equation}\label{Eqn:CMR:Recons:LSB:IP}
\psi_\tau+\vec
W(x,T-\tau)\cdot\nabla\psi=\varepsilon\nabla^2\psi,\quad
\psi(x,0)=\phi(x,T).
\end{equation}
The above PDE mimics thresholding in the sense that it penalizes
the reconstruction from $\vec W$ by introducing a vanishing
viscosity $\varepsilon\nabla^2\psi$, which forces some information
outside the inpainting region flows into the inpainting regions.
Also, when $\varepsilon\to 0$, the solution of
\eqref{Eqn:CMR:Recons:LSB:IP} converges to the viscosity solution
of \eqref{Eqn:CMR:Recons:LSB} \cite{CL1, CIL}.

Since we generally expect volumes of surfaces to increase during
inpainting, we choose the following PDE for the MST,
\begin{equation}\label{Eqn:LSMotion:Inpainting}
\phi_t+(c+\kappa_a-\kappa)|\nabla\phi|=0,\quad
\phi(x,0)=\phi_0(x), \quad c>0.
\end{equation}
Note that the PDE \eqref{Eqn:LSMotion:Inpainting} generates a mean
curvature motion with increasing volumes of the domains enclosed
by level sets of $\phi(x,t)$. The constant $c$ can be regarded as
a parameter that needs to be adjusted according to different
surface inpainting scenarios. In our experiments, we solve PDE
\eqref{Eqn:LSMotion:Inpainting} via a combination of Algorithm
\ref{MotionInNormalDirection_Alg} and Algorithm
\ref{VPMeanCurvMotion_Alg} recalled in Section
\ref{Section:LSBMSR:Numerical}

\begin{algorithm}
\caption{Surface Inpainting via MSR\label{Inpainting_Alg_LSB}}
\begin{algorithmic}
\STATE Start from $\phi_0$, with inpainting region $D$. Choose
some $\varepsilon>0$.

\WHILE{``Not converge"}

\STATE{\textbf{1.}} Perform discrete MST by solving
\eqref{Eqn:LSMotion:Inpainting} and acquire $\vec W_{|i}$ by
Algorithm \ref{MSTIMS LSB_Alg}.

\STATE{\textbf{2.}} Perform IMST by solving
\eqref{Eqn:CMR:Recons:LSB:IP} and obtain $\psi_\varepsilon$.

\STATE{\textbf{3.}} Copy the known information to
$\psi_\varepsilon$: $\psi_\varepsilon |_{D^c}\leftarrow \psi_0
|_{D^c}$.

\STATE{\textbf{4.}} Decrease amount of smoothing:
$\varepsilon\searrow$.

\ENDWHILE
\end{algorithmic}
\end{algorithm}

\begin{figure}[ht]
\centering
    \includegraphics[width=4.5in]{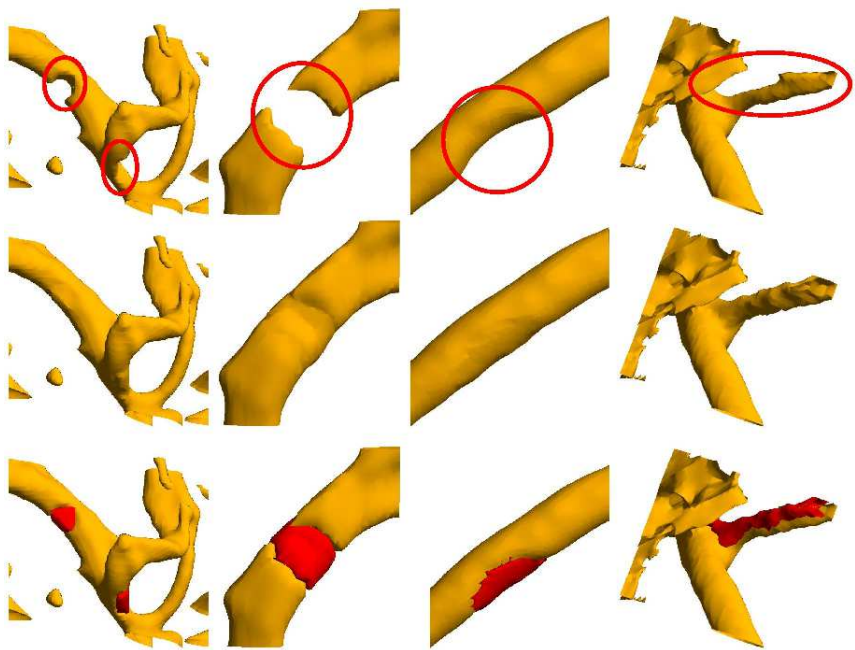}
    \caption{Blood vessel inpainting. Row 1: vessels before inpainting;
    row 2: vessels after inpainting; row 3: inpainted regions shown in red. The percentages of the volume of inpainted region over
    that of the entire shape are: $5.3\%$, $19.2\%$, $6.7\%$ and $5.7\%$.}\label{Fig:Inpainting}
\end{figure}

We test Algorithm \ref{Inpainting_Alg_LSB} on both phantom (first
two vessels in Figure \ref{Fig:Inpainting}) and real (last two
vessels in Figure \ref{Fig:Inpainting}) surface inpainting
scenarios. First row of Figure \ref{Fig:Inpainting} shows four
blood vessels with inpainting regions specified by red circles.
For the two phantom inpainting scenarios, the inpainting regions
are created manually, and the surface within those regions were
chopped off. For the two real inpainting scenarios, we do not know
the exact inpainting regions. Therefore in practice, we adopt a
user interactive strategy to determine the inpainting regions.
After several points have been selected on the surface, the
inpainting regions are then generated automatically. Inpainting
results are given in second and third row of Figure
\ref{Fig:Inpainting}. We want to point out that during the
inpainting process, topological change may occur for some cases
(e.g. second vessel in Figure \ref{Fig:Inpainting}). Although it
violates the assumption in Proposition \ref{Prop:LSBMSR},
topological change is still allowed for inpainting problems. The
reason is that perfect reconstruction is only required at the very
last stage of inpainting (i.e. when $\varepsilon\approx0$) in
order to ensure convergence, while topological changes will happen
during the middle of the process if the parameters (e.g. $c$ in
\eqref{Eqn:LSMotion:Inpainting}) are properly chosen.

\section{Conclusion}

In this paper, we introduced a novel multiscale representation
(MSR) for shapes which is intrinsic to the shape itself, does not
need any parametrization, and the details of the MSR reveals
important geometric information. Based on the MSR, we then
proposed a surface inpainting algorithm and applied it to recover
corrupted blood vessels. This technique is especially useful to
study arteriosclerosis and vessel occlusions. Numerical results
showed that the inpainting regions were nicely filled in according
to the neighboring geometry of the vessels and allowed us to
accurately estimate the volume loss of vessels.

There are still many interesting aspects of both the MSR itself
and its applications worth discovering. For example, a rigorous
analysis of how Algorithm \ref{MSTIMS LSB_Alg} approximates the
continuous version in Definition \ref{Def:CMR:Alt} needs to be
done. Also, we can apply our MSR framework to other type of shape
processing and analysis problems, e.g. shape registration and
classification. In fact we believe that many techniques based on
classical MSR (or multiresolution analysis) for functions (e.g.
wavelets) can now be extended to shapes.

\bibliographystyle{amsplain}

\end{document}